\documentclass[12pt]{article}

\usepackage{graphicx,amssymb,amsmath,amsthm,amsfonts}
\usepackage[top=3cm,bottom=3cm,right=3cm,left=3cm]{geometry}

\begin{document} 

\begin{center}
\footnotesize Communications in Nonlinear Science and Numerical Simulation\\
Volume 69, April 2019, Pages 73-77. DOI: 10.1016/j.cnsns.2018.09.011
\end{center}

\

\begin{center}
\large \bf Lie point symmetries and conservation laws for a class of BBM-KdV systems
\end{center}

\begin{center}
\large Valter Aparecido Silva Junior \!$^{a,b}$
\end{center}

\begin{flushleft}
\footnotesize $^a$ Instituto Federal de Educa\c c\~ao, Ci\^encia e Tecnologia de S\~ao Paulo - IFSP. Avenida Marginal, 585, Fazenda Nossa Senhora Aparecida do Jaguari, 13871-298, S\~ao Jo\~ao da Boa Vista, SP, Brasil.
\end{flushleft}

\begin{flushleft}
\footnotesize $^b$ Instituto de F\'isica ``Gleb Wataghin" - IFGW, Universidade Estadual de Campinas - UNICAMP. Rua S\'ergio Buarque de Holanda, 777, Cidade Universit\'aria, 13083-859, Campinas, SP, Brasil.
\end{flushleft} 

\begin{center}
  \large{\textbf{Abstract}}
\end{center}

We determine the Lie point symmetries of a class of BBM-KdV systems and establish its nonlinear self-adjointness. We then construct conservation laws via Ibragimov's Theorem.

\vspace{.5cm}
\noindent{\textbf{Keywords}: BBM-KdV system, Lie point symmetries, nonlinear self-adjointness, conservation laws.}
\section{Introduction}

Motivated by the works \cite{bcs1, bcs2}, we introduce the following class of systems
\begin{equation}
\label{BBM-KdV}
\left\{
\begin{aligned}
&F_1 \equiv u_t + (a + b)vu_x + (au + c)v_x + \epsilon u_{txx} + \kappa v_{xxx} = 0\\
&F_2 \equiv v_t + (bu + c)u_x + (a + b)vv_x + \lambda u_{xxx} + \sigma v_{txx} = 0
\end{aligned}
\right.,
\end{equation}
henceforth simply referred to as BBM-KdV system, a two-component generalization\footnote{If $u = v$, the system (\ref{BBM-KdV}) is reduced to equations $u_t + [(2a + b)u + c]u_x + \epsilon u_{txx} + \kappa u_{xxx} = 0$ and $u_t + [(a + 2b)u + c]u_x + \sigma u_{txx} + \lambda u_{xxx} = 0$. Both contain (\ref{BBM}) and (\ref{KdV}) as special cases.} of the classic equations \cite{bbm}
\begin{equation}
\label{BBM}
BBM\!\!: \ u_t + (u + 1)u_x - u_{txx} = 0
\end{equation}
and
\begin{equation}
\label{KdV}
KdV\!\!: \ u_t + (u + 1)u_x + u_{xxx} = 0,
\end{equation}
with the objective of studying it from the point of view of the group analysis.

In (\ref{BBM-KdV}), the constants are such that $(a + b)c \neq 0$ and $\{\epsilon, \kappa, \lambda, \sigma\} \neq \{0\}$. Particularly when $a = c = 1$ and $b = 0$, we obtain the already widely investigated systems of Boussinesq ($\epsilon = \kappa = \lambda = 0$, $\sigma = -1/3$), Kaup ($\epsilon = \lambda = \sigma = 0$, $\kappa = 1/3$) and Bona-Smith ($\epsilon = \sigma = \lambda/2 - 1/6$, $\kappa = 0$, $\lambda < 0$), all of them first-order approximations to the Euler equations in the framework of hydrodynamics. Useful in situations where dissipative effects are not significant, these models provide a good description for the two-dimensional motion of small-amplitude long waves on the surface of an ideal fluid. In this context, the independent variable $x$ represents the distance traveled along a fixed depth channel and $t$ the time. The quantities $u(t, x)$ and $v(t, x)$ are related to the deviation of the surface from its undisturbed level and to the horizontal velocity of the fluid, respectively. For more information, see \cite{bcs1,bcs2} and references therein. Relevant results, including exact solutions, can be found in \cite{adm,c,ddlm}.

It's well known that evolution equations don't possess an usual Lagrangian. Therefore this paper is thus organized: first we determine the Lie point symmetries (Section 2) of the BBM-KdV system and establish its nonlinear self-adjointness (Section 3); we then construct conservation laws via Ibragimov's Theorem (Section 4), an extension of the celebrated Noether's Theorem to problems with no variational structure. In the next sections, unless otherwise stated, $c_i$'s are arbitrary constants. All functions are smooth.

We consider that the reader is familiar with the fundamental concepts of group analysis. The basic literature used is \cite{bk,gan,i1,i3,i2,i4,i5,ol}.

\section{Lie Point Symmetries Classification}

Without many details, applying the standard algorithm presented in \cite{bk} and \cite{ol}, a differential operator
\begin{equation}
X = \mathcal{T}(t, x, u, v)\frac{\partial}{\partial t} + \mathcal{X}(t, x, u, v)\frac{\partial}{\partial x} + \mathcal{U}(t, x, u, v)\frac{\partial}{\partial u} + \mathcal{V}(t, x, u, v)\frac{\partial}{\partial v}\nonumber
\end{equation}
generates the Lie point symmetries of the system (\ref{BBM-KdV}) if the conditions of invariance (the so-called determining equations)
\begin{equation}
\label{determining}
\begin{gathered}
\mathcal{T}_x = \mathcal{T}_u = \mathcal{T}_v = \mathcal{X}_u = \mathcal{X}_v = 0,\\
\epsilon\mathcal{X}_t = \epsilon\mathcal{X}_x = \sigma\mathcal{X}_t = \sigma\mathcal{X}_x = 0,\\
\mathcal{U}_t = \mathcal{U}_x = \mathcal{U}_v = \mathcal{V}_t = \mathcal{V}_x = a\mathcal{U} - (au + c)\mathcal{U}_u = 0,\\
b\mathcal{U} + (bu + c)[\mathcal{U}_u + 2(\mathcal{T}_t - \mathcal{X}_x)] = 0,\\
\kappa(\mathcal{U}_u + 2\mathcal{X}_x) = \lambda[\mathcal{U}_u + 2(\mathcal{T}_t - 2\mathcal{X}_x)] = 0,\\
(a + b)[\mathcal{V} + (\mathcal{T}_t - \mathcal{X}_x)v] - \mathcal{X}_t = 0
\end{gathered}
\end{equation}
are satisfied. From (\ref{determining}), it's easy to see that
\begin{equation}
\begin{gathered}
\mathcal{T} = (a + b)c_1t + c_2, \ \mathcal{X} = (a + b)(c_3x + c_4t) + c_5,\\
\mathcal{U} = 2(c_3 - c_1)(au + c), \ \mathcal{V} = (a + b)(c_3 - c_1)v + c_4
\end{gathered}
\nonumber
\end{equation}
with
\begin{equation}
\begin{gathered}
b(a - b)(c_1 - c_3) = 0, \ \epsilon c_3 = \epsilon c_4 = \kappa[ac_1 - (2a + b)c_3] = 0,\\
\lambda[bc_1 - (a + 2b)c_3]= \sigma c_3 = \sigma c_4 = 0.
\end{gathered}
\nonumber
\end{equation}

\

\noindent{{\bf Proposition 1.} \textit{The Lie point symmetries of the BBM-KdV system are summarized in Table 1, where}
\begin{equation}
\begin{gathered}
X_1 = (a + b)\left(t\frac{\partial}{\partial t} - v\frac{\partial}{\partial v}\right) - 2(au + c)\frac{\partial}{\partial u},\\
X_2 = \frac{\partial}{\partial t}, \ X_3 = (a + b)\left(x\frac{\partial}{\partial x} + v\frac{\partial}{\partial v}\right) + 2(au + c)\frac{\partial}{\partial u},\\
X_4 = (a + b)t\frac{\partial}{\partial x} + \frac{\partial}{\partial v}, \ X_5 = \frac{\partial}{\partial x}.
\end{gathered}
\nonumber
\end{equation}

\begin{table}[h!]
\centering
\begin{tabular}{ccccl}
\cline{2-4}
\multicolumn{1}{c|}{}  & \multicolumn{1}{c|}{$b = 0$} & \multicolumn{1}{c|}{$a = b$} & \multicolumn{1}{c|}{$b(a - b) \neq 0$} &  \\ \cline{1-4}
\multicolumn{1}{|c|}{$\{\epsilon, \sigma\} = \{0\}$} & \multicolumn{1}{c|}{\begin{tabular}{clll}
$X_1$ ($\kappa = 0)$ \\
$2X_1 + X_3$ ($\lambda = 0$) \\
$X_2$, $X_4$, $X_5$
\end{tabular}} & \multicolumn{1}{c|}{\begin{tabular}{cllll}
$3X_1 + X_3$ \\
$X_2$, $X_4$, $X_5$
\end{tabular}} & \multicolumn{1}{c|}{$X_2$, $X_4$, $X_5$} &  \\ \cline{1-4}
\multicolumn{1}{|c|}{$\{\epsilon, \sigma\} \neq \{0\}$} & \multicolumn{1}{c|}{\begin{tabular}{cllll}
$X_1$ ($\kappa = 0$), $X_2$, $X_5$
\end{tabular}} & \multicolumn{1}{c|}{\begin{tabular}{cllll}
$X_1$ ($\kappa = \lambda = 0$) \\
$X_2$, $X_5$
\end{tabular}} & \multicolumn{1}{c|}{$X_2$, $X_5$} &  \\ \cline{1-4}
                       & \multicolumn{1}{l}{}  & \multicolumn{1}{l}{}  & \multicolumn{1}{l}{}  &
\end{tabular}
\\
\tablename{ 1.}
\end{table}

\section{Self-Adjointness Classification}

To begin with, let $\bar{u}$ and $\bar{v}$ be the new dependent variables. The formal Lagrangian of the system (\ref{BBM-KdV}) is
\begin{equation}
\mathcal{L} = \bar{u}F_1 + \bar{v}F_2.
\nonumber
\end{equation}

Calculated the adjoint equations
\begin{equation}
\left\{
\begin{aligned}
&F_1^* \equiv -\frac{\delta\mathcal{L}}{\delta u} = \bar{u}_t + (a + b)v\bar{u}_x + \ \! (bu + c)\bar{v}_x + \ \! b\bar{u}v_x + \epsilon\bar{u}_{txx} + \lambda\bar{v}_{xxx} = 0\\
&F_2^* \equiv -\frac{\delta\mathcal{L}}{\delta v} = \bar{v}_t + (au + c)\bar{u}_x + (a + b)v\bar{v}_x - b\bar{u}u_x + \kappa\bar{u}_{xxx} + \sigma\bar{v}_{txx} = 0
\end{aligned}
\right.,\nonumber
\end{equation}
where $\delta/\delta u$ and $\delta/\delta v$ are Euler-Lagrange operators, we assume that
\begin{equation}
F_1^*|_{(\bar{u}, \bar{v}) = (\varphi, \psi)} = MF_1 + NF_2, \quad F_2^*|_{(\bar{u}, \bar{v}) = (\varphi, \psi)} = PF_1 + QF_2.\label{self}
\end{equation}
Here $M$, $N$, $P$ and $Q$ is a set of coefficients to be determined and
\begin{equation}
\label{sub}
\varphi = \varphi(t, x, u, v), \quad \psi = \psi(t, x, u, v)
\end{equation}
two functions that not vanish simultaneously. As
\begin{equation}
F_1^*|_{(\bar{u}, \bar{v}) = (\varphi, \psi)} = D_t\varphi + (a + b)vD_x\varphi + (bu + c)D_x\psi + b\varphi v_x + \epsilon D_tD_x^2\varphi + \lambda D_x^3\psi
\nonumber
\end{equation}
and
\begin{equation}
F_2^*|_{(\bar{u}, \bar{v}) = (\varphi, \psi)} = D_t\psi + (au + c)D_x\varphi + (a + b)vD_x\psi - b\varphi u_x + \kappa D_x^3\varphi + \sigma D_tD_x^2\psi,
\nonumber
\end{equation}
from (\ref{self}) it's possible to conclude that $M = \varphi_u$, $N = \varphi_v$, $P = \psi_u$, $Q = \psi_v$ and
\begin{equation}
\label{self-determining}
\begin{gathered}
\varphi_t + (a + b)v\varphi_x = \epsilon\varphi_x = 0,\\
\psi_t + (au + c)\varphi_x = \psi_x = 0,\\
b\varphi = (au + c)\varphi_u - (bu + c)\psi_v,\\
\varphi_v - \psi_u = (\epsilon - \sigma)\varphi_v = \kappa\varphi_u - \lambda\psi_v = 0,\\
\epsilon\varphi_{uu} = \varphi_{uv} = \varphi_{vv} = 0.
\end{gathered}
\nonumber
\end{equation}

Hence
\begin{equation}
\varphi = (c_1t + c_2)av + f(u) - c_1x, \quad \psi = c_3v + (c_1t + c_2)au + c_1ct + c_4\nonumber
\end{equation}
with
\begin{equation}
\begin{gathered}
bc_1 = bc_2 = \epsilon c_1 = \sigma c_1 = (\epsilon - \sigma)c_2 = 0,\\
\epsilon f''(u) = \kappa f'(u) - \lambda c_3 = 0,\\
bf(u) = (au + c)f'(u) - (bu + c)c_3.
\end{gathered}
\nonumber
\end{equation}

\

\noindent{{\bf Proposition 2.}} \textit{The BBM-KdV system is nonlinearly self-adjoint. The substitutions $(\ref{sub})$ are as follows.}

\

\noindent{\textbf{i)}} \textit{If $b = 0$,}
\begin{equation}
\varphi = (c_1t + c_2)av + c_3c\ln(au + c) - c_1x + c_4, \quad \psi = c_3av + (c_1t + c_2)(au + c) + c_5\nonumber
\end{equation}
\textit{where}
\begin{equation}
\left\{
\begin{aligned}
&c_1 = 0, \ \mathrm{to} \ \{\epsilon, \sigma\} \neq \{0\},\\
&c_2 = 0, \ \mathrm{to} \ \epsilon \neq \sigma,\\
&c_3 = 0, \ \mathrm{to} \ \{\epsilon, \kappa, \lambda\} \neq \{0\}.
\end{aligned}
\right.\nonumber
\end{equation}

\noindent{\textbf{ii)}} \textit{If $a = b$,}
\begin{equation}
\varphi = (au + c)[c_1\ln(au + c) + c_2], \quad \psi = c_1av + c_3\nonumber
\end{equation}
\textit{where}
\begin{equation}
\left\{
\begin{aligned}
&c_1 = 0, \ \mathrm{to} \ \{\epsilon, \kappa, \lambda\} \neq \{0\},\\
&c_2 = 0, \ \mathrm{to} \ \kappa \neq 0.\\
\end{aligned}
\right.\nonumber
\end{equation}

\noindent{\textbf{iii)}} \textit{Let $b(a - b) \neq 0$.}

\textbf{iii.a)} \textit{If $a = 0$,}
\begin{equation}
\varphi = c_1e^{bu/c} - c_2(bu + 2c), \quad \psi = c_2bv + c_3\nonumber
\end{equation}

\textit{where}
\begin{equation}
\left\{
\begin{aligned}
&c_1 = 0, \ \mathrm{to} \ \{\epsilon, \kappa\} \neq \{0\},\\
&c_2 = 0, \ \mathrm{to} \ \kappa \neq -\lambda.
\end{aligned}
\right.\nonumber
\end{equation}

\textbf{iii.b)} \textit{If $a \neq 0$,}
\begin{equation}
\varphi = c_1(au + c)^{b/a} + c_2[b^2u + (2b - a)c], \quad \psi = c_2(a - b)bv + c_3\nonumber
\end{equation}

\textit{where}
\begin{equation}
\left\{
\begin{aligned}
&c_1 = 0, \ \mathrm{to} \ \{\epsilon, \kappa\} \neq \{0\},\\
&c_2 = 0, \ \mathrm{to} \ \lambda a \neq (\kappa + \lambda)b.
\end{aligned}
\right.\nonumber
\end{equation}

\

\noindent{\textbf{Remark.} Actually, the system (\ref{BBM-KdV}) is quasi self-adjoint. It becomes strictly self-adjoint in only two circumstances: $a = 2b$ and $\kappa = \lambda$; or $b = 0$ and $\epsilon = \sigma$.}

\section{Conservation Laws}

In view of Proposition 2, the components of the conserved vector $C = (C^t, C^x)$ associated to $X$, a Lie point symmetry admitted by the system (\ref{BBM-KdV}), are according to Ibragimov's Theorem given by
\begin{equation}
\label{componente-t}
C^t = (\varphi - \epsilon D_x\varphi D_x)W^u + (\psi - \sigma D_x\psi D_x)W^v\nonumber
\end{equation}
and
\begin{equation}
\begin{split}
C^x = [(a &+ b)v\varphi + (bu + c)\psi + \epsilon(\varphi D_tD_x + D_tD_x\varphi) + \lambda(\psi D_x^2 - D_x\psi D_x + D_x^2\psi)]W^u +\\ 
&+ [(au + c)\varphi + (a + b)v\psi + \sigma(\psi D_tD_x + D_tD_x\psi) + \kappa(\varphi D_x^2 - D_x\varphi D_x + D_x^2\varphi)]W^v,
\end{split}
\nonumber
\end{equation}
with
\begin{equation}
W^u = \mathcal{U} - \mathcal{T}u_t - \mathcal{X}u_x, \quad W^v = \mathcal{V} - \mathcal{T}v_t - \mathcal{X}v_x.\nonumber
\end{equation}

We find the conservation laws corresponding to each generator of Table 1. In most cases, however, we are led to trivial vectors or the vectors
\begin{equation}
C^t = u + \epsilon u_{xx}, \quad C^x = (au + c)v + \kappa v_{xx}\nonumber
\end{equation}
and
\begin{equation}
C^t = 2(v + \sigma v_{xx}), \quad C^x = (a + b)v^2 + (bu + 2c)u + 2\lambda u_{xx}\nonumber
\end{equation}
that can be obtained from the first (when b = 0) and second equation of the BBM-KdV system by simple integration (obvious conservation laws). The really interesting cases we list below.

\

\noindent{{\bf Proposition 3.} \textbf{i)} {\it Let $b = 0$.}

\

\textbf{i.a)} {\it From $X_1$, $2X_1 + X_3$ and $X_2$, we obtain}
\begin{equation}
\begin{gathered}
C^t = 2(uv - \epsilon u_xv_x),\\
C^x = cu^2 + (2au + c)v^2 - (\lambda u_x^2 + \kappa v_x^2) + 2[u(\lambda u_{x} + \epsilon v_{t})_x + v(\epsilon u_{t} + \kappa v_{x})_x]
\end{gathered}
\nonumber
\end{equation}

{\it when $\epsilon = \sigma$.}

\

\textbf{i.b)} {\it For $\epsilon = \kappa = 0$, $X_1$ also provides}
\begin{equation}
\begin{gathered}
C^t = \frac{1}{a}(au + c)\ln(au + c) + \frac{a}{2c}(v^2 - \sigma v_x^2),\\
C^x = (au + c)[\ln(au + c) + 1]v + \frac{av}{c}\!\left(\frac{av^2}{3} + \sigma v_{tx}\right)
\end{gathered}
\nonumber
\end{equation}

{\it when $\lambda = 0$ and}
\begin{equation}
\begin{gathered}
C^t = 2[t(au + c)v - xu],\\
C^x = t[c(au + 2c)u - a\lambda u_x^2] + 2(au + c)[(atv - x)v + \lambda tu_{xx}]
\end{gathered}
\nonumber
\end{equation}

{\it when $\sigma = 0$.}

\

\noindent{\textbf{ii)} {\it Let $a = b$.}}

\

\textbf{ii.a)} {\it From $X_1$ and $3X_1 + X_3$, we obtain}
\begin{equation}
\begin{gathered}
C^t = (au + 2c)u - a\epsilon u_x^2,\\
C^x = 2(au + c)[(au + c)v + \epsilon u_{tx}]
\end{gathered}
\nonumber
\end{equation}

{\it when $\kappa = 0$.}

\

\textbf{ii.b)} {\it $X_1$ also provides}
\begin{equation}
\begin{gathered}
C^t = \frac{1}{a}(au + c)^2\ln(au + c) + a(v^2 - \sigma v_x^2),\\
C^x = (au + c)^2[2\ln(au + c) + 1]v + 2av\!\left(\frac{2av^2}{3} + \sigma v_{tx}\right)
\end{gathered}
\nonumber
\end{equation}

{\it when $\epsilon = \kappa = \lambda = 0$.}

\end{document}